\documentclass[11pt]{article}
\usepackage{amsmath}
\usepackage{amssymb}
\usepackage[usenames]{color}
 \usepackage{colortbl}
\def\C{\mathbb{C}}
\def\R{\mathbb{R}}
\def\N{\mathbb{N}}
\def\a{\mathbf{a}}
\def\b{\mathbf{b}}

\newtheorem{theorem}{\hspace*{\parindent}Theorem}
\newtheorem{lemma}{\hspace*{\parindent}Lemma}
\newtheorem{corollary_t}{\hspace*{\parindent}Corollary}[theorem]
\newtheorem{corollary_l}{\hspace*{\parindent}Corollary}[lemma]

\newtheorem{conjecture}{\hspace*{\parindent}Conjecture}
\newcounter{theremark}

\title{Inequalities for series in $q$-shifted factorials and $q$-gamma functions}
\author{S.I.\:Kalmykov$^{\rm a,c}$ and D.B.\:Karp$^{\rm b,c}$\footnote{Corresponding author. E-mail: S.I.\:Kalmykov -- \emph{sergeykalmykov@inbox.ru}, D.B.\:Karp -- \emph{dimkrp@gmail.com}}
\\[10pt]\small{\textit{$\phantom{1}^a$Shanghai Jiao Tong University, Shanghai, China}}\\\small{\textit{$\phantom{1}^b$Far Eastern Federal University, Vladivostok, Russia}}
\\\small{\textit{$\phantom{1}^c$Institute of Applied Mathematics, FEBRAS, Vladivostok, Russia}}}
\date{}
\begin{document}
\maketitle

\begin{abstract}
The paper studies logarithmic convexity and concavity of power series with coefficients involving $q$-gamma functions or $q$-shifted factorials with respect to a parameter contained in their arguments.  The principal motivating examples of such series are basic hypergeometric functions.  We consider four types of series.  For each type we establish conditions sufficient for the power series coefficients of the generalized Tur\'{a}nian formed by these series to have constant sign. Finally, we furnish seven examples of basic hypergeometric functions satisfying our general theorems. This investigation extends our previous results on power series with coefficient involving the ordinary gamma functions and the shifted factorials.
\end{abstract}

\bigskip


\bigskip


\section{Introduction}
In a series of papers \cite{KK1,KK2,KK3,KarpPOMI2011,KarpPOMI2012,KarpPOMI2015,KS}  we considered the logarithmic convexity and concavity with respect to the parameter $\mu$ for the class of functions  defined by the series
\begin{equation}\label{eq:f-general}
f(\mu;x)=A_0(\mu)\sum\limits_{n\geq0}f_nA_1(n+\mu)x^n,
\end{equation}
where the coefficients $A_0(\cdot)$, $A_1(\cdot)$ are chosen from the following nomenclature
\begin{equation}\label{eq:nomenclature}
A_0,A_1\in\left\{1,\Gamma(\cdot),\frac{1}{\Gamma(\cdot)}, \frac{\Gamma(a+\cdot)}{\Gamma(b+\cdot)}\right\}
\end{equation}
and $f_n$ is a (usually non-negative) real sequence.  Here $\Gamma$ stands for Euler's gamma function and $a$, $b$ are non-negative parameters.  The main motivating examples of functions of the form (\ref{eq:f-general}) are the (generalized) hypergeometric functions. Moreover, their  derivatives with respect to parameters other than $\mu$ are also instances of  (\ref{eq:f-general}). Further examples of (\ref{eq:f-general}) can be given  with $f_n$ not expressible in terms of gamma functions.  Our papers \cite{KK1,KK2,KK3,KarpPOMI2011,KarpPOMI2012,KarpPOMI2015,KS} cover nearly all possible combinations of $A_0$ and $A_1$ from the collection (\ref{eq:nomenclature}). Most our results are shaped as follows.  Logarithmic concavity (convexity) of $\mu\to{f(\mu;x)}$ on an interval $I$ is equivalent to non-negativity (non-positivity) of the generalized Tur\'{a}nian
\begin{equation}\label{eq:delta-defined}
\Delta_f(\alpha,\beta;x)=f(\mu+\alpha;x)f(\mu+\beta;x)-f(\mu;x)f(\mu+\alpha+\beta;x)=\sum\limits_{m=0}^{\infty}\delta_mx^m
\end{equation}
for arbitrary $\alpha,\beta\ge0$ such that $\mu,\mu+\alpha,\mu+\beta,\mu+\alpha+\beta\in{I}$.  In many cases, however, we were only able to prove non-negativity of $\Delta_f(\alpha,\beta;x)$ for $\mu,\beta\ge0$ and $\alpha\in\N$, so that our results in such cases are incomplete in the sense that they can probably be still extended to all $\alpha\ge0$.  On the other hand, in most cases we actually demonstrate that, under certain restrictions, all coefficients $\delta_m$ have the same sign.  This type of results can be termed ''coefficient-wise logarithmic concavity/convexity'' and can be viewed as a strengthening of usual log-concavity/log-convexity.

The purpose of this paper is to extend our previous results by substituting the nomenclature (\ref{eq:nomenclature}) with $\{1,\Gamma_q(\cdot)$,  $[\Gamma_q(\cdot)]^{-1}\}$, where $\Gamma_q(\cdot)$ denotes the $q$-gamma function, defined in (\ref{eq:q-Gamma-def}) below.  To this end, we prove four theorems corresponding to four nontrivial  combinations  of $A_0$, $A_1$ chosen from the above set.  We also present a number of corollaries giving two-sided bounds and integral representations for the generalized Tur\'{a}nians (\ref{eq:delta-defined}) and certain product ratios of functions (\ref{eq:f-general}).  Finally, we furnish seven examples of $q$-hypergeometric functions satisfying our general theorems.  Some results dealing with the Tur\'{a}n type inequalities for $q$-hypergeometric functions have been recently obtained by Baricz,  Raghavendar and Swaminathan in \cite{BRS,BS} and Mehrez and Sitnik in \cite{Mehrez,MehrSitn}. In particular, our Theorem~1 can be seen as a far-reaching generalization of \cite[Theorem~3.2]{BRS}, their connection explored in Example~2. Furthermore, our Theorem~3 generalizes some statements from \cite[Theorem~3.1]{BRS} and \cite[Theorem~1]{MehrSitn} which we explore in Example~3.  Continued fractions for and the mapping properties of the ratios of the basic hypergeometric functions have been recently studied in \cite{AS,BS}.

\section{Definitions and preliminaries}

Let us fix some notation and terminology. We will use the standard symbols $\N$, $\R$ and $\C$ to denote natural, real and complex numbers, respectively; $\N_0=\N\cup\{0\}$, $\R_{+}=[0,\infty)$.  A positive function is called logarithmically concave (convex) if its logarithm is concave (convex). Next, a function $f:I\to(0,\infty)$ defined on an interval $I\subset(0,\infty)$ is said to be multiplicatively concave if
$$
f(x^{\lambda}y^{1-\lambda})\ge f(x)^{\lambda}f(y)^{1-\lambda}
$$
for $\lambda\in[0,1]$ and all $x$, $y$ such that $x^{\lambda}y^{1-\lambda}\in{I}$. It is multiplicatively convex if the above inequality is reversed.  In other words this says that $\log(f)$ is concave function of $\log(x)$, i.e. $\log[f(e^x)]$ is concave.  If $f$ is continuous, its multiplicative concavity is equivalent to
\begin{equation}\label{eq:GGconcavity}
f(\sqrt{xy})\ge\sqrt{f(x)f(y)},~~~x,y\in{I},
\end{equation}
which can be termed Jensen multiplicative concavity, GG-concavity or concavity with respect to geometric means \cite[section~2.3]{NP}.
We will need the following elementary lemma.
\begin{lemma}\label{lm:multconcave}
A positive function $f$ is Jensen multiplicatively concave \emph{(}convex\emph{)} on an interval $I\subset(0,\infty)$ iff
\begin{equation}\label{eq:qGGconcavity}
f(a)^2\ge (\le) f(a/q)f(aq)~~\text{for all}~q\in(0,1)~\text{and all}~a~\text{such that}~aq,aq^{-1}\in{I}.
\end{equation}
\end{lemma}
\textbf{Proof.}  Suppose (\ref{eq:GGconcavity}) holds.  By symmetry we may assume that $0<y<x$.
Denote $a=\sqrt{xy}\in{I}$, $q=\sqrt{y/x}\in(0,1)$. Solving these equations we get $y=aq$, $x=a/q$, so that (\ref{eq:GGconcavity}) squared becomes (\ref{eq:qGGconcavity}).  For reverse implication note that $aq,aq^{-1}\in{I}$ implies $a\in{I}$ and apply the inverse change of variable.  $\hfill\square$

From the above lemma we conclude that for continuous $f$ (and this is the only case we are dealing with here) inequality (\ref{eq:qGGconcavity}) is equivalent to multiplicative concavity.  Recall that a nonnegative function $f$ defined on an interval $I$ is called completely monotonic there if it has derivatives of all orders and $(-1)^nf^{(n)}(x)\ge0$ for $n\in\N_0$ and $x\in{I}$, see \cite[Defintion~1.3]{SSV}.

\begin{lemma}\label{lm:multconvex}
Suppose $\phi(x)=\sum_{k\ge0}\phi_kx^k$ converges for $|x|<R$ with $0<R\le\infty$ and $\phi_k\ge0$. Then $x\to\phi(x)$ is multiplicatively convex and $y\to\phi(1/y)$ is completely monotonic on $(1/R,\infty)$.
\end{lemma}
\textbf{Proof.} Hardy, Littlewood and P\'{o}lya theorem \cite[Proposition 2.3.3]{NP}
states that functions with non-negative power series coefficients are multiplicatively convex.
Furthermore, $y^{-m}$, $m\in\N_0$, is apparently completely monotonic and a convergent series of completely monotonic functions with nonnegative coefficients is again completely monotonic according to \cite[Theorem 3]{MS}. $\hfill\square$

A sequence $f:\N_0\to\R_+$ is ${\rm PF_2}$ (P\'{o}lya frequency sub two) or doubly positive if it is nontrivial, log-concave, $f_k^2\geq f_{k-1}f_{k+1}$,
$k\in\N$, and has no internal zeros. The last claim means that $f_N=0$ implies either $f_{N+i}=0$ for all $i\in \mathbb{N}_0$ or $f_{N-i}=0$ for $i=0,\ldots,N$.

The next two lemmas can be found in \cite[Lemmas~2 and 3]{KK2}.
\begin{lemma}\label{lm:alphaone}
Let $f$ be a nonnegative-valued function defined on $\R_{+}$ and
$$
\Delta_{f}(\alpha,\beta)=f(\mu+\alpha)f(\mu+\beta)-f(\mu)f(\mu+\beta+\alpha)\ge0~ \text{for}~\alpha=1~\text{and all}~\mu,\beta\ge0.
$$
Then $\Delta_{f}(\alpha,\beta)\ge0$ for all $\alpha\in\N$ and $\mu,\beta\ge0$. If inequality is strict in the hypothesis of the lemma then it is also strict in the conclusion.
\end{lemma}

\begin{lemma}\label{lm:nncoef}
Let $f=\sum_{n=0}^{\infty}f_n(\mu)x^{n}$ and $\Delta_f(\alpha,\beta;x)$ be defined in \emph{(\ref{eq:delta-defined})}. Suppose that $\Delta_f(1,\beta;x)$ has non-negative power series coefficients for all $\mu,\beta\ge0$. Then $\Delta_f(\alpha,\beta;x)$ has non-negative power series coefficients for all $\alpha\in\N$, $\alpha\le\beta+1$ and $\mu,\beta\ge0$. If the coefficients in the hypothesis are strictly positive, then they are strictly positive in the conclusion.
\end{lemma}

Next, we formulate an elementary inequality we will repeatedly use below.
\begin{lemma}\label{lm:uvrs}
Suppose $u,v,r,s>0$, $u=\max(u,v,r,s)$ and $uv>rs$. Then
$u+v>r+s$.
\end{lemma}
The proof is straightforward and will be omitted. Note also that Lemma~\ref{lm:uvrs} is a particular case of a much more general result on logarithmic majorization, see \cite[2.A.b]{MOA}.

In the next lemma proved in \cite[Lemma~2.1]{KK1} we say that a sequence has no more than one change of sign if it has the pattern $(--\cdots--00\cdots00++\cdots++)$, where zeros and minus signs may be missing.
\begin{lemma}\label{lm:sum}
Suppose $\{f_k\}_{k=0}^{n}$ is a doubly positive sequence and $A_0,A_1,\ldots,A_{[n/2]}$ is a real sequence satisfying $A_{[n/2]}>0$,
$\sum\limits_{0\leq{k}\leq{n/2}}\!\!\!\!A_k\geq{0}$ and having no more than one change of sign. Then
\begin{equation*}\label{eq:keysum}
\sum\limits_{0\leq{k}\leq{n/2}}f_{k}f_{n-k}A_k\geq{0}.
\end{equation*}
Equality is only attained if $f_k=f_0\alpha^k$, $\alpha>0$, and
$\sum\limits_{0\leq{k}\leq{n/2}}\!\!\!\!A_k=0$.
\end{lemma}

We will use the standard definition of the $q$-shifted factorial \cite{GR,KacCheung}:
$$
(a;q)_0=1,~~~~(a;q)_n=\prod\limits_{k=0}^{n-1}(1-aq^k),~~n\in\N.
$$
This definition works for any complex $a$ and $q$ but in this paper we confine ourselves to the case $0<q<1$.
Under this restriction we also define
$$
(a;q)_{\infty}=\lim\limits_{n\to\infty}(a;q)_n,
$$
where the limit can be shown to exist as a finite number for all complex $a$. For products of the $q$-shifted factorials we will use the usual
short-hand notation
$$
(a_1, a_2,\dots, a_m; q)_{n}=(a_1;q)_{n}(a_2;q)_{n}\cdots (a_m;q)_{n},
$$
where $n$ may take the value $\infty$.  The $q$-gamma function is defined by
\begin{equation}\label{eq:q-Gamma-def}
\Gamma_q(z)=(1-q)^{1-z}\frac{(q;q)_{\infty}}{(q^z;q)_{\infty}}
\end{equation}
for $|q|<1$ and all complex $z$ such that $q^{z+k}\ne1$ for $k\in\N_0$. A $q$-number $[x]_q$ is the ratio $(1-q^x)/(1-q)$. The generalized $q$-hypergeometric series is defined by \cite[formula (1.2.22)]{GR}
\begin{equation}\label{eq:r-phi-s}
{_r\phi_s}(a_1,a_2,\dots, a_r;b_1,\dots,b_s;q;z)=\sum_{n=0}^{\infty}\frac{(a_1,a_2,\dots,a_r;q)_n}{(q,b_1,\dots,b_s;q)_n} \left[(-1)^nq^{\binom{n}{2}}\right]^{1+s-r}z^n,
\end{equation}
where $r\le{s+1}$ and the series converges for all $z$ if $r\le{s}$ and for $|z|<1$ if $r=s+1$ \cite[section~1.2]{GR}.

In our previous work  on log-concavity we repeatedly used the celebrated Chu-Vandermonde identity \cite[formula (1.2.9)]{GR} as, for instance, in \cite[Theorem~1]{KS} and \cite[Lemma~2.2]{KK1}.   This identity can be obtained by equating coefficients in $(1-z)^a(1-z)^b=(1-z)^{a+b}$.  Here, we will need a $q$-analogue of this simple equality.
The $q$-analogue usually found in the literature reads \cite[formula (1.3.13)]{GR}
$$
{_1\phi_0}(a;-;z){_1\phi_0}(b;-;az)={_1\phi_0}(ab;-;z).
$$
The apparent asymmetry of the left-hand side with respect to $a$ and $b$ makes this formula inappropriate for our purposes. Instead, we will use the following symmetric version.
\begin{lemma}\label{lm:qbinomial}
The following $q$-identity holds when both sides are well defined
\begin{equation}\label{eq:qbinomial}
{_1\phi_0}(a;-;z){_1\phi_0}(b;-;z)={_1\phi_0}(ab;-;z){_2\phi_1}(a,b;abz;z).
\end{equation}
\end{lemma}
\textbf{Proof.} Applying $q$-binomial theorem \cite[formula (1.3.2)]{GR} and Heine's $q$-Gauss theorem
\cite[formula (1.5.1)]{GR} we get:
\begin{multline*}
{_1\phi_0}(a;-;z){_1\phi_0}(b;-;z)=\frac{(az,bz;q)_{\infty}}{[(z;q)_{\infty}]^2}
\\
=\frac{(abz;q)_{\infty}}{(z;q)_{\infty}}
\frac{(az,bz;q)_{\infty}}{(abz,z;q)_{\infty}}={_1\phi_0}(ab;-;z){_2\phi_1}(a,b;abz;z).~~\square
\end{multline*}

The next lemma is a $q$-analogue of \cite[Lemma~2.3]{KK1}.
\begin{lemma}\label{lm:rec-q-gammas}
Suppose $m\geq{0}$ is an integer. Then for all complex $\beta$ and
$\mu$
\begin{multline}\label{eq:rec-q-gamma-sum}
S_m(\mu,\beta)=\sum\limits_{k=0}^{m}\left\{\frac{1}{\Gamma_q(k+\mu+1)\Gamma_q(m-k+\mu+\beta)}
-\frac{1}{\Gamma_q(k+\mu)\Gamma_q(m-k+\mu+\beta+1)}\right\}
\\
=\frac{(q^{\mu+\beta};q)_{m+1}-(q^\mu;q)_{m+1}}{\Gamma_q(\mu+m+1)\Gamma_q(\mu+\beta+m+1)(1-q)^{m+1}}.
\end{multline}
\end{lemma}
\textbf{Proof.}   We will use the following simple properties of $q$-gamma function which are straightforward from its definition  (\ref{eq:q-Gamma-def}):
\begin{equation}\label{eq:qGamma-prop}
\Gamma_q(x+1)=[x]_{q}\Gamma_q(x),~~~~~~\frac{\Gamma_q(x+k)}{\Gamma_q(x)}=\frac{(q^x;q)_k}{(1-q)^k}.
\end{equation}
Denote $a=q^\mu$, $b=q^{\beta}$ for brevity and compute utilizing the above formulas,
\begin{multline*}
S_m(\mu,\beta)=\frac{1}{\Gamma_q(\mu+1)\Gamma_q(\mu+\beta+1)}\\
\times\sum\limits_{k=0}^{m}\left\{\frac{\Gamma_q(\mu+1)\Gamma_q(\mu+\beta)[\mu+\beta]_q}{\Gamma_q(\mu+1+k)\Gamma_q(\mu+\beta+m-k)}
-\frac{\Gamma_q(\mu)\Gamma_q(\mu+\beta+1)[\mu]_{q}}{\Gamma_q(\mu+k)\Gamma_q(\mu+\beta+1+m-k)}\right\}
\\
=\frac{(1-q)^m}{\Gamma_q(\mu+1)\Gamma_q(\mu+\beta+1)}
\sum\limits_{k=0}^{m}\left\{\frac{[\mu+\beta]_q}{(q^{\mu+1};q)_k(q^{\mu+\beta};q)_{m-k}}
-\frac{[\mu]_{q}}{(q^\mu;q)_k (q^{\mu+\beta+1};q)_{m-k}}\right\}
\\
=\frac{(1-q)^{m-1}}{\Gamma_q(\mu+1)\Gamma_q(\mu+\beta+1)}
\sum\limits_{k=0}^{m}\frac{(1-a)(1-ab)}{(a;q)_k(ab;q)_{m-k}}\left\{\frac{1}{1-aq^k}
-\frac{1}{1-abq^{m-k}}\right\}
\\
=\frac{a(1-q)^{m-1}}{\Gamma_q(\mu+1)\Gamma_q(\mu+\beta+1)}
\sum\limits_{k=0}^{m}\frac{q^k-bq^{m-k}}{(aq;q)_k(abq;q)_{m-k}}.
\end{multline*}
Define
$$
u_k=\frac{1}{(aq;q)_{k-1}(abq;q)_{m-k}},~~1\leq{k}\leq{m},
~~u_0=\frac{1-a}{(abq;q)_{m}},~~u_{m+1}=\frac{1-ab}{(aq;q)_{m}}.
$$
An easy calculation shows that
$$
u_{k+1}-u_{k}=\frac{a(q^k-bq^{m-k})}{(aq;q)_{k}(abq;q)_{m-k}}
$$
for $k=0,1,\ldots,m$, so that
$$
a\sum\limits_{k=0}^{m}\frac{q^k-bq^{m-k}}{(aq;q)_k(abq;q)_{m-k}}=\sum\limits_{k=0}^{m}(u_{k+1}-u_{k})=u_{m+1}-u_{0}
=\frac{(ab;q)_{m+1}-(a;q)_{m+1}}{(aq;q)_{m}(abq;q)_{m}}
$$
and
$$
S_m(\mu,\beta)=\frac{(1-q)^{m-1}}{\Gamma_q(\mu+1)\Gamma_q(\mu+\beta+1)}\frac{(ab;q)_{m+1}-(a;q)_{m+1}}{(aq;q)_{m}(abq;q)_{m}},
$$
which is equivalent to (\ref{eq:rec-q-gamma-sum}) after substituting back $a=q^\mu$, $b=q^{\beta}$.~$\hfill\square$

\smallskip

The next corollary is a straightforward consequence of formula (\ref{eq:rec-q-gamma-sum}).
\begin{corollary_l}\label{cr:rec-q-Gamma}
If $m\in\N_0$, $\mu\geq{-1}$, $\beta\geq{0}$ and $\mu+\beta\geq{0}$, then $S_{m}(\mu,\beta)\geq{0}$.  The inequality is strict unless $\beta=0$.
\end{corollary_l}

\section{Main results}

Our first theorem is a $q$-analogue of \cite[Theorem~1]{KS}.  The power series in this theorem is to be understood as formal.
Nonetheless, we will show in a remark below that  it has a positive radius of convergence under the hypotheses of the theorem.
\begin{theorem}\label{th:pochhammertop}
For a real sequence $\{f_n\}_{n\ge0}$ and fixed $0<q<1$ define
\begin{equation}\label{eq:f}
f(a;x)=\sum\limits_{n=0}^{\infty}f_n\frac{(a;q)_n}{(q;q)_n}x^n.
\end{equation}
Suppose $\{f_n\}_{n\ge0}$ is doubly positive. Then the generalized Tur\'{a}nian
\begin{equation}\label{eq:phiab}
\Delta_{f}(\alpha,\beta;x)=f(q^{\mu+\alpha};x)f(q^{\mu+\beta};x)-f(q^{\mu};x)f(q^{\mu+\alpha+\beta};x)
\end{equation}
has positive coefficients at all positive powers of $x$ for all $\alpha,\beta>0$ and $\mu\ge0$.  In particular, $\mu\to{f(q^{\mu};x)}$
is log-concave on $[0,\infty)$ for each fixed $x>0$ in the domain of convergence.
\end{theorem}
\textbf{Proof.} Writing $\Delta_{f}(\alpha,\beta;x)=\sum_{m\geq0}\delta_f(m)x^m$ we calculate using the Cauchy product:
$$
\delta_f(m)=\sum\limits_{k=0}^{m}f_kf_{m-k}\frac{(q^{\mu+\alpha};q)_k(q^{\mu+\beta};q)_{m-k}-(q^{\mu+\alpha+\beta};q)_k(q^{\mu};q)_{m-k}}{(q;q)_k(q;q)_{m-k}}
=\sum\limits_{k=0}^{m}f_kf_{m-k}M_k,
$$
where the last equality is the definition of the numbers $M_k$.  We aim at applying Lemma~\ref{lm:sum} with $A_k=M_k+M_{m-k}$ for $0\le{k}<m/2$, $A_{k}=M_{k}$ for $k=m/2$ (this term is only present for even $m$). To this end, we first  need to demonstrate that
$$
\widehat{\delta}_f(m)=\sum\limits_{0\leq{k}\leq{m/2}}\!\!\!\!A_k=\sum\limits_{k=0}^{m}M_k>0.
$$
This definition together with (\ref{eq:r-phi-s}) and Lemma~\ref{lm:qbinomial} yields the next chain of equalities:
\begin{multline*}
\sum\limits_{m=0}^{\infty}\widehat{\delta}_f(m)x^m={_{1}\phi_0}(q^{\mu+\alpha};-;x) {_{1}\phi_0}(q^{\mu+\beta};-;x)-{_{1}\phi_0}(q^{\mu+\alpha+\beta};-;x){_{1}\phi_0}(q^{\mu};-;x)
\\
={_{1}\phi_0}(q^{2\mu+\alpha+\beta};-;x)\left[{_{2}\phi_1(q^{\mu+\alpha},q^{\mu+\beta};q^{2\mu+\alpha+\beta}x;x)}-{_{2}\phi_1(q^{\mu+\alpha+\beta},q^{\mu};q^{2\mu+\alpha+\beta}x;x)}\right]
\\
={_{1}\phi_0}(q^{2\mu+\alpha+\beta};-;x)\sum_{k=0}^{\infty}\frac{(q^{\mu+\alpha},q^{\mu+\beta};q)_k-(q^{\mu+\alpha+\beta},q^{\mu};q)_k}{(q^{2\mu+\alpha+\beta}x;q)_k(q;q)_k}x^k.
\end{multline*}
As
$$
\frac{1}{(q^{2\mu+\alpha+\beta}x;q)_k}=\frac{(q^{2\mu+\alpha+\beta+k}x;q)_{\infty}}{(q^{2\mu+\alpha+\beta}x;q)_{\infty}}=\sum\limits_{n=0}^{\infty}\frac{(q^{2\mu+\alpha+\beta+k};q)_{n}}{(q;q)_{n}}x^n
$$
by the $q$-binomial theorem \cite[formula (1.3.2)]{GR}, to prove positivity of $\widehat{\delta}_f(m)$ it remains to check that the difference $(q^{\mu+\alpha},q^{\mu+\beta};q)_k-(q^{\mu+\alpha+\beta},q^{\mu};q)_k$ is positive. This amounts to
\begin{equation*}
\begin{split}
(q^{\mu+\alpha},q^{\mu+\beta};q)_k-(q^{\mu+\alpha+\beta},q^{\mu};q)_k
&=\frac{(q^{\mu+\alpha},q^{\mu+\beta};q)_{\infty}}{(q^{\mu+\alpha+k},q^{\mu+\beta+k};q)_{\infty}}
-\frac{(q^{\mu},q^{\mu+\alpha+\beta};q)_{\infty}}{(q^{\mu+k},q^{\mu+\alpha+\beta+k};q)_{\infty}}>0
\\
\Leftarrow~~~&\frac{(1-q^{\alpha}s_j)(1-q^{\beta}s_j)}{(1-s_j)(1-q^{\alpha+\beta}s_j)}>\frac{(1-q^{\alpha}t_j)(1-q^{\beta}t_j)}{(1-t_j)(1-q^{\alpha+\beta}t_j)},
\end{split}
\end{equation*}
where $s_j=q^{\mu+j}>t_j=q^{\mu+k+j}$ for all $j\in\N_0$ and $k\in\N$.  The last inequality is true since the function
\begin{equation}\label{eq:U}
x\mapsto U(x)=\frac{(1-q^{\alpha}x)(1-q^{\beta}x)}{(1-q^{\alpha+\beta}x)(1-x)}
\end{equation}
is easily seen to be increasing on $[0,1)$. This completes the proof of positivity of  $\widehat{\delta}_f(m)$ for all $m\ge1$.

Our next goal is to show that the sequence $A_0,A_1,\ldots,A_{[m/2]}$ satisfies $A_{[m/2]}>0$ and has no more than one change of sign. It suffices to prove the implication $A_k\leq{0}~\Rightarrow~A_{k-1}<0$ for $k\ge1$.  Indeed, as $\widehat{\delta}_f(m)=\sum_{0\leq{k}\leq{m/2}}A_k>0$, this implication immediately leads to the conclusion that $A_{[m/2]}>0$.  Next,
we spell out:
\begin{multline*}
(q;q)_k(q;q)_{m-k}A_k=\underbrace{(q^{\mu+\alpha};q)_k(q^{\mu+\beta};q)_{m-k}}_{u_k}+\underbrace{(q^{\mu+\alpha};q)_{m-k}(q^{\mu+\beta};q)_{k}}_{v_k}
\\
-\underbrace{(q^{\mu+\alpha+\beta};q)_k(q^{\mu};q)_{m-k}}_{r_k}-\underbrace{(q^{\mu+\alpha+\beta};q)_{m-k}(q^{\mu};q)_{k}}_{s_k}~\text{for}~0\le{k}<m-k.
\end{multline*}
Assuming that $A_k=u_k+v_k-r_k-s_k\le0$, we need to show that
\begin{multline*}
0>(q;q)_{k-1}(q;q)_{m-k+1}A_{k-1}\!=\!\underbrace{\frac{1-q^{\mu+\beta+m-k}}{1-q^{\mu+\alpha+k-1}}}_{=I_1}u_k+\underbrace{\frac{1-q^{\mu+\alpha+m-k}}{1-q^{\mu+\beta+k-1}}}_{=I_2}v_k\\
-\underbrace{\frac{1-q^{\mu+m-k}}{1-q^{\mu+\alpha+\beta +k-1}}}_{=I_3}r_k
-\underbrace{\frac{1-q^{\mu+\alpha+\beta+m-k}}{1-q^{\mu+k-1}}}_{=I_4}s_k.
\end{multline*}
Performing some elementary calculations and employing the increase of the function $U(x)$ from (\ref{eq:U}), we see that $r_k<\min(u_k,v_k,s_k)$, $u_kv_k>r_ks_k$ and $u_k\ge{v_k}$ if $\alpha\le\beta$ while $u_k<v_k$ if $\alpha>\beta$. Combined with our hypothesis $u_k+v_k-r_k-s_k\leq0$ these inequalities imply that $s_k>u_k\ge{v_k}>r_k$ if $\alpha\le\beta$ or $s_k>v_k>{u_k}>r_k$ if $\alpha>\beta$. Indeed, if $u_k\ge{s_k}\geq{v_k}>r_k$  or $u_k\ge{v_k}\geq{s_k}>r_k$, then $u_k+v_k-r_k-s_k>0$.  Rearranging the above expression for $A_{k-1}$ gives
\begin{multline}\label{Iuvrs}
(q;q)_{k-1}(q;q)_{m-k+1}A_{k-1}=
\\
I_2(u_{k}+v_{k}-r_{k}-s_{k})+ (I_2-I_4)(s_k-u_{k})+(I_1-I_4)(u_{k}-r_{k})+(I_1+I_2-I_3-I_4)r_{k}.
\end{multline}
The first term on the right is non-positive by the assumption $A_{k}\leq0$. We will show that all
further terms on the right hand side are negative.  The second term is negative since $s_k>u_k$ (as we have just proved) and
$$
I_2<I_4~\Leftrightarrow~\frac{1-q^{\mu+\alpha+m-k}}{1-q^{\mu+\beta+k-1}}<\frac{1-q^{\mu+\alpha+\beta+m-k}}{1-q^{\mu+k-1}},
$$
which is true as $\beta>0$ and $0<q<1$.  Next,
$$
I_1<I_4~\Leftrightarrow~\frac{1-q^{\mu+\beta+m-k}}{1-q^{\mu+\alpha+k-1}}<\frac{1-q^{\mu+\alpha+\beta+m-k}}{1-q^{\mu+k-1}},
$$
which is true as $\alpha>0$ and $0<q<1$. In view of $u_k>r_k$ the third term in (\ref{Iuvrs}) is then negative.
It remains to show that $I_1+I_2<I_3+I_4$ which will be accomplished by means of Lemma~~\ref{lm:uvrs}.  We have
$$
I_3<I_4~\Leftrightarrow~\frac{1-q^{\mu+m-k}}{1-q^{\mu+\alpha+\beta+k-1}}<\frac{1-q^{\mu+\alpha+\beta+m-k}}{1-q^{\mu+k-1}}
$$
as $\alpha+\beta>0$ and $0<q<1$.  Therefore, $I_4=\max\{I_1,I_2,I_3,I_4\}$.  Furthermore,
$$
I_1 I_2<I_3 I_4 \Leftrightarrow \frac{1-q^{\mu+\alpha+m-k}}{1-q^{\mu+\beta+k-1}}\cdot\frac{1-q^{\mu+\beta+m-k}}{1-q^{\mu+\alpha+k-1}}<\frac{1-q^{\mu+m-k}}{1-q^{\mu+\alpha+\beta+k-1}}
\cdot\frac{1-q^{\mu+\alpha+\beta+m-k}}{1-q^{\mu+k-1}}
$$
$$
\Leftrightarrow~\frac{1-q^{\mu+\alpha+m-k}}{1-q^{\mu+m-k}}\cdot\frac{1-q^{\mu+\beta+m-k}}{1-q^{\mu+\alpha+\beta+m-k}}<\frac{1-q^{\mu+\alpha+k-1}}{1-q^{\mu+k-1}}
\cdot\frac{1-q^{\mu+\beta+k-1}}{1-q^{\mu+\alpha+\beta+k-1}}.
$$
The last inequality holds because  $0\le{k-1}<m-k$ and the function $U(x)$ from (\ref{eq:U}) is increasing on $[0,1)$.  Hence, $I_1+I_2<I_3+I_4$ by Lemma~\ref{lm:uvrs} and negativity of the last term on the right hand side of (\ref{Iuvrs}) follows.  This completes the proof that $A_0,A_1,\ldots,A_{[m/2]}$ has no more than one change of sign and $A_{[m/2]}>0$. According to Lemma~\ref{lm:sum} the coefficients $\delta_f(m)$ are positive for all $m\ge1$.  $\hfill\square$

\smallskip

\textbf{Remark.}  Log-concavity of the sequence $\{f_n\}_{n\ge0}$ implies that the series (\ref{eq:f}) has a positive radius of convergence, $R_f>0$. Indeed, by log-concavity we have $f_{n}/f_{n+1}\ge{f_{n-1}/f_{n}}$, so that the nonnegative sequence  $\{f_{n}/f_{n+1}\}_{n\ge0}$ is increasing. Hence,
$$
R_f=\lim_{n\to\infty}\left|\frac{f_n(a;q)_n(q;q)_{n+1}}{f_{n+1}(a;q)_{n+1}(q;q)_n}\right|=\lim_{n\to\infty}\left|\frac{f_n(1-q^{n+1})}{f_{n+1}(1-aq^{n})}\right|=
\lim_{n\to\infty}\left|\frac{f_n}{f_{n+1}}\right|>0.
$$

\smallskip

\textbf{Remark.} The conclusion of Theorem~\ref{th:pochhammertop} holds trivially also for $-2k-1\le\mu<-2k$, $k\in\N_0$, as long as $\mu+\alpha$ and $\mu+\beta$ are positive, since in this case the rightmost term in (\ref{eq:phiab}) has non-positive power series coefficients.

\begin{corollary_t}\label{cr:Delta-fCM}
Suppose $\{f_n\}_{n\ge0}$ is a doubly positive sequence. Then the function $a\to{f(a;x)}$ defined in \emph{(\ref{eq:f})} is multiplicatively concave on $(0,1)$. If also $\alpha,\beta>0$,  and $\mu\ge0$, then the function $x\to\Delta_{f}(\alpha,\beta;x)$ defined in \emph{(\ref{eq:phiab})} is multiplicatively convex on $(0,R_f)$, where $R_f$ is the radius of convergence in \emph{(\ref{eq:f})}, while the function $y\to\Delta_f(\alpha,\beta;1/y)$ is completely monotonic \emph{(}and therefore log-convex\emph{)} on $(1/R_f,\infty)$.
\end{corollary_t}
\textbf{Proof.} Multiplicative concavity of $a\to{f(a;x)}$ follows from Theorem~\ref{th:pochhammertop} by Lemma~\ref{lm:multconcave}. Multiplicative convexity of
 $x\to\Delta_{f}(\alpha,\beta;x)$ and complete monotonicity of $y\to\Delta_f(\alpha,\beta;1/y)$ are implied by Lemma~\ref{lm:multconvex} as the power series coefficients of $\Delta_f(\alpha,\beta;x)$ are positive by Theorem~\ref{th:pochhammertop}. $\hfill\square$

\smallskip

\textbf{Remark}.   Complete monotonicity of $y\to\Delta_f(\alpha,\beta;1/y)$ implies by Bernstein's theorem \cite[Theorem 1.4]{SSV} that there exists a non-negative measure $\tau$ supported on $[0,\infty)$ such that
$$
\Delta_f(\alpha,\beta;x)=\int_{[0,\infty)}e^{-(1/x-1/R_f)t}\tau(dt).
$$
If $R_f=\infty$, this measure is given by
$$
\tau(dt)=\delta_f(0)\mathbf{1}_{0}+\left(\sum\nolimits_{m=1}^{\infty}\frac{\delta_f(m)t^m}{(m-1)!}\right)dt,
$$
where $\mathbf{1}_{0}$ is the unit mass concentrated at zero. This formula can be easily verified by termwise integration.  Note, that in this situation the function
$\Delta_f(\alpha,\beta;1/y)$ satisfies the conditions of \cite[Theorem~1.1]{KoumPeders} and hence enjoys all the properties stated in that theorem.

\smallskip

The next theorem is a $q$-analogue of \cite[Theorem~2]{KS}.  Again, the power series in this theorem is to be understood as formal.

\begin{theorem}\label{th:q-gamma-num}
Suppose $\{d_n\}_{n\ge0}$ is a non-negative sequence, $0<q<1$ is fixed and $d(\mu;x)$ is defined by
\begin{equation}\label{eq:d-def}
d(\mu;x)=\sum\limits_{n=0}^{\infty}d_n\Gamma_q(\mu+n)x^n.
\end{equation}
Then, the generalized Tur\'{a}nian
$$
\Delta_d(\alpha,\beta;x)=d(\mu+\alpha;x)d(\mu+\beta;x)-d(\mu;x)d(\mu+\alpha+\beta;x)
$$
has negative coefficients at all powers of $x$ for all $\mu,\alpha,\beta>0$.  In particular, the function $\mu\to{d(\mu;x)}$ is log-convex on $(0,\infty)$ for each fixed $x>0$ in the domain of convergence.
\end{theorem}
\textbf{Proof.} Writing $\Delta_d(\alpha,\beta;x)=\sum_{m=0}^{\infty}\delta_d(m)x^m$ we calculate using the Cauchy product:
$$
\delta_d(m)=\sum\limits_{k=0}^{m}d_kd_{m-k}\{\Gamma_q(\mu+\alpha+k)\Gamma_q(\mu+\beta+m-k)-\Gamma_q(\mu+k)\Gamma_q(\mu+\alpha+\beta+m-k)\}.
$$
We can rewrite $\delta_d(m)$ in the form
$$
\delta_d(m)=\sum\limits_{k=0}^{[m/2]}d_kd_{m-k}A_k,
$$
where, for $k<m/2$,
\begin{multline*}
A_k=\underbrace{\Gamma_q(\mu+\alpha+k)\Gamma_q(\mu+\beta+m-k)}_{=u_k}
+\underbrace{\Gamma_q(\mu+\alpha+m-k)\Gamma_q(\mu+\beta+k)}_{=v_k}
\\
-\underbrace{\Gamma_q(\mu+k)\Gamma_q(\mu+\alpha+\beta+m-k)}_{=r_k}
-\underbrace{\Gamma_q(\mu+m-k)\Gamma_q(\mu+\alpha+\beta+k)}_{=s_k}
\end{multline*}
and, for $k=m/2$,
$$
A_k=\Gamma_q(\mu+\alpha+k)\Gamma_q(\mu+\beta+m-k)-\Gamma_q(\mu+k)\Gamma_q(\mu+\alpha+\beta+m-k).
$$
We aim to show that $A_k<0$ for $k=0,\ldots,[m/2]$ using Lemma~\ref{lm:uvrs}.  For $k<m-k$ the following comparisons between
the numbers $u_k$, $v_k$, $r_k$ and $s_k$ are straightforward to verify from the increase of $x\to\Gamma_q(\mu+\gamma+x)/\Gamma_q(\mu+x)$ for any $\gamma>0$ (which is equivalent to log-convexity of $x\to\Gamma_q(x)$):
$$
u_k<r_k~\Leftrightarrow~\frac{\Gamma_q(\mu+\alpha+k)}{\Gamma_q(\mu+k)}<\frac{\Gamma_q(\mu+\alpha+\beta+m-k)}{\Gamma_q(\mu+\beta+m-k)},
$$
$$
v_k<r_k~\Leftrightarrow~\frac{\Gamma_q(\mu+\beta+k)}{\Gamma_q(\mu+k)}<\frac{\Gamma_q(\mu+\alpha+\beta+m-k)}{\Gamma_q(\mu+\alpha+m-k)},
$$
$$
s_k<r_k\Leftrightarrow~\frac{\Gamma_q(\mu+\alpha+\beta+k)}{\Gamma_q(\mu+k)}<\frac{\Gamma_q(\mu+\alpha+\beta+m-k)}{\Gamma_q(\mu+m-k)}
$$
and
$$
u_kv_k<r_ks_k~\Leftrightarrow~\frac{\Gamma_q(\mu+\beta+k)}{\Gamma_q(\mu+k)}\cdot\frac{\Gamma_q(\mu+\beta+m-k)}{\Gamma_q(\mu+m-k)}\!<\! \frac{\Gamma_q(\mu+\alpha+\beta+k)}{\Gamma_q(\mu+\alpha+k)}\cdot\frac{\Gamma_q(\mu+\alpha+\beta+m-k)}{\Gamma_q(\mu+\alpha+m-k)}.
$$
Therefore, according to Lemma~\ref{lm:uvrs} $u_k+v_k<r_k+s_k$, so that  $A_k<0$ for $0<k<m/2$. For even $m$ the inequality between $u_k$ and $r_k$ remains true
and implies $A_k<0$ for $k=m/2$. $\hfill\square$

\smallskip
In the next two corollaries $R_d>0$ denotes the radius of convergence of the series in (\ref{eq:d-def}).  Our first corollary is
similar to Corollary~\ref{cr:Delta-fCM} and is a direct consequence of Lemma~\ref{lm:multconvex}.

\begin{corollary_t}\label{cr:Delta-dCM}
Under conditions of Theorem~\ref{th:q-gamma-num} and for all $\alpha,\beta,\mu>0$,  the function $x\to-\Delta_d(\alpha,\beta;x)$  is multiplicatively convex on $(0,R_d)$, while the function $y\to-\Delta_d(\alpha,\beta;1/y)$ is completely monotonic \emph{(}and therefore log-convex\emph{)} on $(1/R_d,\infty)$.
\end{corollary_t}

\smallskip

The next two corollaries are obtained by joint application of Theorems~\ref{th:pochhammertop} and \ref{th:q-gamma-num}.

\begin{corollary_t}\label{cr:Delta-d-bounds}
Suppose $\{d_n(q;q)_n\}_{n\ge0}$ is a doubly positive sequence and $d(\mu;x)$ is defined in \emph{(\ref{eq:d-def})}.  Then for all $\mu,\alpha,\beta>0$ and $0\le{x}<R_d$
the following estimates hold\emph{:}
$$
\frac{\Gamma_q(\mu+\alpha)\Gamma_q(\mu+\beta)}{\Gamma_q(\mu)\Gamma_q(\mu+\alpha+\beta)}\le\frac{d(\mu+\alpha;x)d(\mu+\beta;x)}{d(\mu;x)d(\mu+\alpha+\beta;x)}<1
$$
and
\begin{multline*}
\left[\frac{\Gamma_q(\mu+\alpha)\Gamma_q(\mu+\beta)}{\Gamma_q(\mu)\Gamma_q(\mu+\alpha+\beta)}-1\right]d(\mu;x)d(\mu+\alpha+\beta;x)\le\Delta_d(\alpha,\beta;x)
\\
\le d_0^2\left[\Gamma_q(\mu+\alpha)\Gamma_q(\mu+\beta)-\Gamma_q(\mu)\Gamma_q(\mu+\alpha+\beta)\right].
\end{multline*}
with equality only at $x=0$.  The upper bounds in both inequalities remain valid if $\{d_n\}_{n\ge0}$ is any non-negative sequence.
\end{corollary_t}
\textbf{Proof.}  The upper bound in the first inequality is equivalent to  $\Delta_d(\alpha,\beta;x)<0$ which is one of the conclusions of Theorem~\ref{th:q-gamma-num}.
To prove the lower bound define $f_n=d_n(q;q)_n(1-q)^{-n}$ and notice that $\{f_n\}_{n\ge0}$ is doubly positive by hypothesis of the corollary. If $f(a;x)$ is given by (\ref{eq:f}) then (\ref{eq:qGamma-prop}) implies that $f(q^{\mu};x)=d(\mu;x)/\Gamma_q(\mu)$. Hence, the lower bound is equivalent to $\Delta_f(\alpha,\beta;x)\ge0$ which is a conclusion of Theorem~\ref{th:pochhammertop}.  The right hand side of the second inequality is equal to $\delta_d(0)$, so that the upper bound in the second inequality follows immediately from Theorem~\ref{th:q-gamma-num}. The lower bound is a rearrangement of the lower bound in the first inequality. $\hfill\square$

\smallskip

Recall that $R_f$ is the radius of convergence in (\ref{eq:f}).  We have
\begin{corollary_t}\label{cr:Delta-f-bounds}
Suppose $\{f_n\}_{n\ge0}$ is a doubly positive sequence and $f(\mu;x)$ is defined in \emph{(\ref{eq:f})}.  Then for all $\alpha,\beta>0$, $\mu\ge0$ and $0\le{x}<R_f$
the following bounds hold true\emph{:}
\begin{equation*}\label{eq:twosided2}
\frac{\Gamma_q(\mu+\alpha)\Gamma_q(\mu+\beta)}{\Gamma_q(\mu)\Gamma_q(\mu+\alpha+\beta)}
<\frac{f(q^{\mu};x)f(q^{\mu+\alpha+\beta};x)}{f(q^{\mu+\alpha};x)f(q^{\mu+\beta};x)}\le1
\end{equation*}
and
$$
f_0f_1(1-q^{\alpha})(1-q^{\beta})\frac{q^{\mu}x}{1-q}\!\le\!\Delta_f(\alpha,\beta;x)\!<\!\left[1-\frac{\Gamma_q(\mu+\alpha)\Gamma_q(\mu+\beta)}{\Gamma_q(\mu)\Gamma_q(\mu+\alpha+\beta)}\right]
\!f(q^{\mu+\alpha};x)f(q^{\mu+\beta};x)
$$
with equality only at $x=0$.  The lower bound in the first inequality and the upper bound in the second inequality remain valid if $\{f_n\}_{n\ge0}$ is any non-negative sequence.
\end{corollary_t}
\textbf{Proof.} The upper bound in the first inequality is equivalent to  $\Delta_f(\alpha,\beta;x)\ge0$ which is one of the conclusions of Theorem~\ref{th:pochhammertop}.
To prove the lower bound define $d_n=(1-q)^nf_n$. If $d(\mu;x)$ is given by (\ref{eq:d-def}), then (\ref{eq:qGamma-prop}) implies that $d(\mu;x)=\Gamma_q(\mu)f(q^{\mu};x)$.
Hence, the lower bound is equivalent to $\Delta_d(\alpha,\beta;x)<0$ which is a conclusion of Theorem~\ref{th:q-gamma-num}.  The left hand side of the second inequality is equal to $\delta_f(1)x$, so that the lower bound follows immediately from Theorem~\ref{th:pochhammertop}. The upper bound is a rearrangement of the lower bound in the first inequality. $\hfill\square$

\smallskip

The next theorem is a $q$-analogue of \cite[Theorem~3.5]{KK1}.
\begin{theorem}\label{th:gammadenom1}
For a real sequence $\{g_n\}_{n\ge0}$ and fixed $0<q<1$ define
\begin{equation}\label{eq:g-def}
g(\mu;x)=\sum\limits_{n=0}^{\infty}\frac{g_nx^n}{\Gamma_q(\mu+n)}.
\end{equation}
Suppose $\{g_n\}_{n=0}^{\infty}$ is doubly positive. Then $(\ref{eq:g-def})$ has a positive radius of convergence $R_g$ and the generalized Tur\'{a}nian
\begin{equation*}\label{eq:tur}
\Delta_g(\alpha,\beta;x)=g(\mu+\alpha;x)g(\mu+\beta;x)-g(\mu;x)g(\mu+\alpha+\beta;x)
\end{equation*}
is positive for $0<x<R_g$ if $\alpha\in\N$, $\beta>0$ and $\mu\ge\max(-\beta,-1)$.  If, in addition, $\alpha\le\beta+1$, then
$\Delta_g(\alpha,\beta;x)$ has positive power series coefficients at all powers of $x$.
\end{theorem}
\textbf{Proof.} The convergence of the series (\ref{eq:g-def}) in some disk when $\{g_n\}_{n=0}^{\infty}$ is doubly positive follows from the argument given in the remark after Theorem~\ref{th:pochhammertop} together with the asymptotic formula for the $q$-gamma function \cite[(2.4)]{Daalhuis}. To prove the remaining statements it is sufficient to consider the case $\alpha=1$ according to Lemmas~\ref{lm:alphaone} and \ref{lm:nncoef}. Writing $\Delta_g(\alpha,\beta;x)=\sum_{m=0}^{\infty}\delta_g(m)x^m$ we calculate using the Cauchy product:
\begin{multline*}
\delta_g(m)=\sum\limits_{k=0}^{m}\left\{\frac{g_kg_{m-k}}{\Gamma_q(k+\mu+1)\Gamma_q(m-k+\mu+\beta)}-\frac{g_kg_{m-k}}{\Gamma_q(m-k+\mu)\Gamma_q(k+\mu+1+\beta)}\right\}
\\
=\sum\limits_{0\leq{k}\leq{m/2}}g_kg_{m-k}A_k,
\end{multline*}
where
\begin{multline*}
A_k=[\Gamma_q(k+\mu+1)\Gamma_q(m-k+\mu+\beta)]^{-1}
+[\Gamma_q(k+\mu+\beta)\Gamma_q(m-k+\mu+1)]^{-1}
\\[0pt]
-[\Gamma_q(m-k+\mu)\Gamma_q(k+\mu+1+\beta)]^{-1}
-[\Gamma_q(k+\mu)\Gamma_q(m-k+\mu+1+\beta)]^{-1}
\end{multline*}
for $k<m/2$ and
$$
A_{m/2}=[\Gamma_q(m/2+\mu+1)\Gamma_q(m/2+\mu+\beta)]^{-1}-[\Gamma_q(m/2+\mu)\Gamma_q(m/2+\mu+1+\beta)]^{-1}
$$
for $k=m/2$ (this term is only present for even $m$).  This last term is always positive since $x\to\Gamma_q(x+\alpha)/\Gamma_q(x)$ is increasing on $(0,\infty)$ for any $\alpha>0$ due to log-convexity of $x\to\Gamma_q(x)$. For $S_m(\mu,\beta)$ defined in (\ref{eq:rec-q-gamma-sum}) we get
\begin{equation*}\label{eq:pos}
\sum\limits_{0\leq{k}\leq{m/2}}A_k=S_m(\mu,\beta)>0
\end{equation*}
by Corollary~\ref{cr:rec-q-Gamma}.  We will demonstrate that the sequence $\{A_k\}_{k=0}^{[m/2]}$ has no more than one change of sign in order to apply Lemma~\ref{lm:sum}. This amounts to showing the implication $A_{k}\leq0~\Rightarrow~A_{k-1}<0$ for $1\leq{k}<m/2$.
Using the second formula in (\ref{eq:qGamma-prop}) we can rewrite $A_k$ as
$$
A_k=\frac{(1-q)^{m+1}}{\Gamma_q(\mu)\Gamma_q(\mu+\beta)}\left(F(a,b,k)+F(a,b,m-k)-F(b,a,k)-F(b,a,m-k)\right),
$$
where $a:=q^{\mu}>b:=q^{\mu+\beta}$ and $F(a,b,k):=[(a;q)_{k+1}(b;q)_{m-k}]^{-1}$.  It is easy to verify that
$$
F(b,a,m-k)=\frac{1-aq^k}{1-bq^{m-k}}F(a,b,k)~~~\text{and}~~~F(a,b,m-k)=\frac{1-bq^k}{1-aq^{m-k}}F(b,a,k),
$$
which yields
$$
A_k=\frac{(1-q)^{m+1}q^k}{\Gamma_q(\mu)\Gamma_q(\mu+\beta)}\left(\frac{a-bq^{m-2k}}{1-bq^{m-k}}F(a,b,k)+\frac{aq^{m-2k}-b}{1-aq^{m-k}}F(b,a,k)\right).
$$
Hence, the condition $A_k\leq0$ is equivalent to
$$
\frac{(b;q)_{k+1}(a;q)_{m-k+1}}{(a;q)_{k+1}(b;q)_{m-k+1}}\leq\frac{b-aq^{m-2k}}{a-bq^{m-2k}}.
$$
Assuming this to be true for some $1\leq{k}<m/2$, we immediately conclude that $A_{k-1}<0$ as
$$
\frac{(b;q)_{k}(a;q)_{m-k+2}}{(a;q)_{k}(b;q)_{m-k+2}}<\frac{(b;q)_{k+1}(a;q)_{m-k+1}}{(a;q)_{k+1}(b;q)_{m-k+1}}\leq\frac{b-aq^{m-2k}}{a-bq^{m-2k}}<\frac{b-aq^{m-2k+2}}{a-bq^{m-2k+2}},
$$
where both the rightmost and the leftmost inequalities follow from $0<b<a<1$. An application of Lemma~\ref{lm:sum} completes the proof.  $\hfill\square$

\smallskip

We believe that the restrictions $\alpha\in\N$ and $\alpha\le\beta+1$ in the hypotheses of Theorem~\ref{th:gammadenom1} solely reflect the limitations of our method of proof. Writing $R_g>0$ for the radius of convergence in (\ref{eq:g-def}), we claim that the next conjecture is true.

\begin{conjecture}\label{con:gammadenom}
All conclusions of Theorem~\ref{th:gammadenom1} are valid for all $\alpha>0$, while conditions on other parameters remain intact.  In particular, $\mu\to{g(\mu;x)}$
is log-concave for each fixed $0<x<R_g$.
\end{conjecture}

\begin{corollary_t}\label{cr:g-multconc}
Suppose $\{g_n\}_{n\ge0}$ is a doubly positive sequence. Then the function
$$
a\to\hat{g}(a;y)=(a;q)_{\infty}\sum\limits_{n=0}^{\infty}\frac{g_ny^n}{(a;q)_n}
$$
is multiplicatively concave on $(0,1)$  for each fixed $0<y<R_g$, where $R_g$ is the radius of convergence of this series.
\end{corollary_t}
\textbf{Proof.}  Indeed, on applying the definition (\ref{eq:q-Gamma-def}) and the obvious identity $(q^{\mu+n};q)_{\infty}=(q^{\mu};q)_{\infty}/(q^{\mu};q)_{n}$,
the inequality $\Delta_g(1,1;x)\ge0$ takes the form
$$
(q^{\mu'};q)_{\infty}^2\left(\sum\limits_{n=0}^{\infty}\frac{g_ny^n}{(q^{\mu'};q)_n}\right)^2\ge(q^{\mu'-1};q)_{\infty}(q^{\mu'+1};q)_{\infty}
\left(\sum\limits_{n=0}^{\infty}\frac{g_ny^n}{(q^{\mu'-1};q)_n}\right)\left(\sum\limits_{n=0}^{\infty}\frac{g_ny^n}{(q^{\mu'+1};q)_n}\right),
$$
where $\mu'=\mu+1$ and $y=x(1-q)$.  Setting $a=q^{\mu'}$ this inequality can be rewritten as $\hat{g}(a;y)^2\geq\hat{g}(a/q;y)\hat{g}(aq;y)$.  The claim now follows by Lemma~\ref{lm:multconcave}. $\hfill\square$

\begin{corollary_t}\label{cr:cm2}
Suppose $\beta>0$, $\beta+1\ge\alpha\in\N$,  $\mu\ge0$. Then the function $x\to\Delta_{g}(\alpha,\beta;x)$  is multiplicatively convex on $(0,R_g)$ and the function $y\to\Delta_g(\alpha,\beta;1/y)$ is completely monotonic \emph{(}and therefore log-convex\emph{)} on $(1/R_g,\infty)$.
\end{corollary_t}

\smallskip

The next theorem is a $q$-analogue of \cite[Theorem~3]{KS}.  The power series in this theorem is to be understood as formal.
It might converge or diverge depending on the behavior of the coefficients.
\begin{theorem}\label{th:q-poch-denom}
Suppose $\{h_n\}_{n\ge0}$ is a non-negative sequence, $0<q<1$ is fixed and $h(a;x)$ is defined by
\begin{equation}\label{eq:h}
h(a;x)=\sum\limits_{n=0}^{\infty}\frac{h_{n}x^n}{(a;q)_n}.
\end{equation}
Then, the generalized Tur\'{a}nian
\begin{equation}\label{eq:Delta-h}
\Delta_h(\alpha,\beta;x)=h(q^{\mu+\alpha};x)h(q^{\mu+\beta};x)-h(q^{\mu};x)h(q^{\mu+\alpha+\beta};x)
\end{equation}
has negative coefficients at all positive powers of $x$ for all $\mu,\alpha,\beta>0$.  In particular, the function $\mu\to{h(q^{\mu};x)}$ is log-convex on $(0,\infty)$ for each fixed $x>0$ in the domain of convergence.
\end{theorem}
\textbf{Proof.} Writing $\Delta_h(\alpha,\beta;x)=\sum_{m=0}^{\infty}\delta_h(m)x^m$ we calculate using the Cauchy product:
\begin{multline*}
-\delta_h(m)=\sum\limits_{k=0}^{m}h_kh_{m-k}\left(\frac{1}{(q^{\mu};q)_k(q^{\mu+\alpha+\beta};q)_{m-k}}-\frac{1}{(q^{\mu+\alpha};q)_k(q^{\mu+\beta};q)_{m-k}}\right)\\
=\sum\limits_{0\leq{k}\leq{m/2}}h_kh_{m-k}A_k,
\end{multline*}
where, for $k<m/2$,
\begin{multline*}
A_k=\underbrace{\frac{1}{(q^{\mu};q)_k(q^{\mu+\alpha+\beta};q)_{m-k}}}_{=u_k}
+\underbrace{\frac{1}{(q^{\mu};q)_{m-k}(q^{\mu+\alpha+\beta};q)_{k}}}_{=v_k}
\\
-\underbrace{\frac{1}{(q^{\mu+\alpha};q)_k(q^{\mu+\beta};q)_{m-k}}}_{=r_k}
-\underbrace{\frac{1}{(q^{\mu+\alpha};q)_{m-k}(q^{\mu+\beta};q)_{k}}}_{=s_k},
\end{multline*}
and, for $k=m/2$,
$$
A_k=\frac{1}{(q^{\mu};q)_k(q^{\mu+\alpha+\beta};q)_{m-k}}-\frac{1}{(q^{\mu+\alpha};q)_k(q^{\mu+\beta};q)_{m-k}}.
$$
The following comparisons between the numbers $u_k$, $v_k$, $r_k$ and $s_k$ are straightforward to verify:
\begin{equation*}
v_k>u_k \Leftrightarrow \frac{1}{(q^{\mu};q)_{m-k}(q^{\mu+\alpha+\beta};q)_{k}}>\frac{1}{(q^{\mu};q)_k(q^{\mu+\alpha+\beta};q)_{m-k}}~\Leftrightarrow~\prod\limits_{i=k}^{m-k-1}\frac{1-q^{\mu+i}}{1- q^{\mu+\alpha+\beta+i}}<1.
\end{equation*}
Next,
\begin{multline*}
v_k>r_k\Leftrightarrow \frac{1}{(q^{\mu};q)_{m-k}(q^{\mu+\alpha+\beta};q)_{k}}> \frac{1}{(q^{\mu+\beta};q)_{m-k}(q^{\mu+\alpha};q)_{k}}\\[7pt]
\Leftrightarrow \prod\limits_{i=0}^{k-1}\frac{1- q^{\mu+\alpha+i}}{1- q^{\mu+\alpha+\beta+i}}>
\prod\limits_{i=0}^{k-1}\frac{1- q^{\mu+i}}{1- q^{\mu+\beta+i}}\prod\limits_{i=k}^{m-k-1}\frac{1- q^{\mu+i}}{1- q^{\mu+\beta+i}}.
\end{multline*}
The rightmost product is clearly less than one. Comparing $i$-th terms in the remaining two products amounts to
$(1- q^{\mu+\alpha+i})(1- q^{\mu+\beta+i})>(1-q^{\mu+\alpha+\beta+i})(1-q^{\mu+i})$ which is equivalent to $(1-q^{\alpha})(1-q^{\beta})>0$.  The inequality
$v_k>s_k$ is proved similarly by exchanging the roles of $\alpha$ and $\beta$. Finally,
\begin{equation*}
u_kv_k>r_ks_k~\Leftrightarrow~\frac{(q^{\mu+\beta};q)_k(q^{\mu+\beta};q)_{m-k}}{(q^{\mu};q)_k(q^{\mu};q)_{m-k}}>\frac{(q^{\mu+\alpha+\beta};q)_k(q^{\mu+\alpha+\beta};q)_{m-k}}{(q^{\mu+\alpha};q)_k(q^{\mu+\alpha};q)_{m-k}}.
\end{equation*}
The last inequality holds because each function
$$
x\to\dfrac{1-xq^{\mu+\beta+i}}{1-xq^{\mu+i}}, ~~~i=0,1,\ldots,
$$
is increasing on $[0,1]$. Now we are in the position to apply Lemma~\ref{lm:uvrs} to conclude that $u_k+v_k>r_k+s_k$ or $A_k>0$, $0<k<m/2$.
For even $m$ from the inequality between $v_k$ and $r_k$ yields $A_k>0$ for $k=m/2$. $\hfill\square$

\smallskip

By Lemma~\ref{lm:multconcave} we see that $a\to{h(a;x)}$ is multiplicatively convex on $(0,1)$ for each fixed $0<x<R_h$, where $R_h$ is  the radius of convergence in (\ref{eq:h}).  Next, we have:\

\begin{corollary_t}
Suppose $\alpha,\beta,\mu>0$. Then the function $x\to-\Delta_h(\alpha,\beta;x)$  is multiplicatively convex on $(0,R_h)$, while the function $y\to-\Delta_h(\alpha,\beta;1/y)$ is completely monotonic \emph{(}and therefore log-convex\emph{)} on $(1/R_h,\infty)$.
\end{corollary_t}

\smallskip

\begin{corollary_t}\label{cr:Delta-h-bounds}
Suppose $\{h_n\}_{n\ge0}$ is a doubly positive sequence and $h(a;x)$ is defined in \emph{(\ref{eq:h})}.  Then for all $\mu,\beta>0$, $\alpha\in\N$ and $0\le{x}<R_h$
the following estimates hold\emph{:}
$$
\frac{\Gamma_q(\mu+\alpha)\Gamma_q(\mu+\beta)}{\Gamma_q(\mu)\Gamma_q(\mu+\alpha+\beta)}<\frac{h(q^{\mu+\alpha};x)h(q^{\mu+\beta};x)}{h(q^{\mu};x)h(q^{\mu+\alpha+\beta};x)}\le1
$$
and
\begin{multline*}
\left[\frac{\Gamma_q(\mu+\alpha)\Gamma_q(\mu+\beta)}{\Gamma_q(\mu)\Gamma_q(\mu+\alpha+\beta)}-1\right]h(q^{\mu};x)h(q^{\mu+\alpha+\beta};x)<\Delta_h(\alpha,\beta;x)
\\
\le\frac{-h_0h_1xq^{\mu}(1-q^{\alpha})(1-q^{\beta})}{(1-q^{\mu})(1-q^{\mu+\alpha+\beta})(1-q^{\mu+\alpha})(1-q^{\mu+\beta})}
\end{multline*}
with equality only at $x=0$.  The upper bounds in both inequalities remain valid if $\{h_n\}_{n\ge0}$ is any non-negative sequence and $\alpha$ is any positive number.
\end{corollary_t}
\textbf{Proof.}  The upper bound in the first inequality is equivalent to  $\Delta_h(\alpha,\beta;x)\le0$ which is one of the conclusions of Theorem~\ref{th:q-poch-denom}.
To prove the lower bound define $g_n=h_n(1-q)^{-n}$ and notice that $\{g_n\}_{n\ge0}$ is doubly positive. If $g(\mu;x)$ is given by (\ref{eq:g-def}) then (\ref{eq:q-Gamma-def}) shows that $g(\mu;x)=h(q^{\mu};x)/\Gamma_q(\mu)$. Hence, the lower bound is equivalent to $\Delta_g(\alpha,\beta;x)>0$ which is a conclusion of Theorem~\ref{th:gammadenom1}.  The right hand side of the second inequality is equal to $\delta_h(1)x$, so that the upper bound in the second inequality follows immediately from Theorem~\ref{th:q-poch-denom}. The lower bound is a rearrangement of the lower bound in the first inequality. $\hfill\square$

\smallskip

\begin{corollary_t}\label{cr:Delta-g-bounds}
Suppose $\alpha\in\N$, $\beta>0$, $\mu\ge\max(-\beta,-1)$, $\{g_n\}_{n\ge0}$ is a doubly positive sequence and $g(\mu;x)$
is defined in \emph{(\ref{eq:g-def})}.  Then the estimates
\begin{equation}\label{eq:twosided1}
\frac{\Gamma_q(\mu+\alpha)\Gamma_q(\mu+\beta)}{\Gamma_q(\mu)\Gamma_q(\mu+\alpha+\beta)}
\le\frac{g(\mu;x)g(\mu+\alpha+\beta;x)}{g(\mu+\alpha;x)g(\mu+\beta;x)}<1
\end{equation}
and
$$
\frac{g_0^2(1-q)^{\alpha}[(q^{\mu};q)_{\alpha}^{-1}-(q^{\mu+\beta};q)_{\alpha}^{-1}]}{\Gamma_{q}(\mu)\Gamma_{q}(\mu+\beta)}
\!\le\Delta_g(\alpha,\beta;x)\le\!\left[1-\frac{(q^{\mu};q)_{\alpha}}{(q^{\mu+\beta};q)_{\alpha}}\right]\!
g(\mu\!+\!\alpha;x)g(\mu\!+\!\beta;x)
$$
hold for $0\le{x}<R_g$, where $R_g$ is the radius of convergence in \emph{(\ref{eq:g-def})}, and equality is only attained at $x=0$.  The lower bound in the second inequality
requires the additional restriction $\alpha\le\beta+1$.
\end{corollary_t}
\textbf{Proof.}  The upper bound in the first inequality is a restatement of Theorem~\ref{th:gammadenom1} since it is equivalent to $\Delta_g(\alpha,\beta;x)>0$. To prove the lower bound we first assume that $\mu>0$ and define $h_n=g_n(1-q)^n$. For $h(a;x)$ given in (\ref{eq:h}) formula (\ref{eq:q-Gamma-def}) shows that $h(q^{\mu};x)=\Gamma_q(\mu)g(\mu;x)$.  The lower bound is then equivalent to $\Delta_h(\alpha,\beta;x)\le0$.  If $\mu=0$, then the left hand side is zero, while the ratio in the middle is nonnegative.
For $\mu=-1$ the lower bound is again zero if $\alpha,\beta\ne1$.  If $\alpha=1$ or $\beta=1$, then the lower bound is finite and negative while for $\alpha=\beta=1$ it is negative and infinite. It remains to consider $-1<\mu<0$.  Since $\Gamma_{q}(\mu)<0$, the lower bound in (\ref{eq:twosided1}) reduces to
$$
\Gamma_q(\mu+\alpha)g(\mu+\alpha;x)\Gamma_q(\mu+\beta)g(\mu+\beta;x)-\Gamma_q(\mu)g(\mu;x)\Gamma_q(\mu+\alpha+\beta)g(\mu+\alpha+\beta;x)>0.
$$
This inequality holds since $\Gamma_q(\mu)g(\mu;x)=\sum_{n=0}^{\infty}(1-q)^ng_nx^n[(q^{\mu};q)_n]^{-1}$ and the $m$-th power series coefficient of the expression on the left hand side equals
$$
(1-q)^m\sum_{k=0}^mg_kg_{m-k}\left\{\frac{1}{(q^{\mu+\alpha};q)_k(q^{\mu+\beta};q)_{m-k}}-\frac{1}{(q^{\mu};q)_k(q^{\mu+\alpha+\beta};q)_{m-k}}\right\}.
$$
Each term of this sum is positive because $(q^{\mu};q)_k<0$ for $k=1,2,\dots,m$, and $(q^{\mu+\beta};q)_m<(q^{\mu+\alpha+\beta};q)_m$ for $k=0$.

To prove the lower bound in the second inequality we employ the second formula in (\ref{eq:qGamma-prop}) to check that it is equal to $\delta_g(0)$. Hence, the lower bound follows immediately from Theorem~~\ref{th:gammadenom1}.   The upper bound is a rearrangement of the lower bound in the first inequality in view of the second formula in (\ref{eq:qGamma-prop}).$\hfill\square$

\section{Applications}

In this section we explore the application of the general theorems from the previous section to derive inequalities for concrete special functions.
Expectedly, the most natural examples come from the basic hypergeometric functions although their $q$-derivatives in parameters or ''mixed'' classical-basic series could also be considered.

\smallskip

\textbf{Example~1.}  The $q$-analogues of the modified Bessel function $I_{\nu}$ were introduced by Ismail \cite[(2.5)]{Ism} and rediscovered by Olshanetskii and Rogov \cite[section~3.1]{OR1}. These $q$-analogues are given by \cite[(17), (18)]{OR1}
$$
I^{(1)}_{\nu}(y;q)=\frac{(y/2)^{\nu}}{(1-q)^{\nu}}\sum_{n=0}^{\infty}\frac{(y/2)^{2n}}{(1-q)^n(q;q)_n\Gamma_q(\nu+n+1)},~~~|y|<2,
$$
and
$$
I^{(2)}_{\nu}(y;q)=\frac{(y/2)^{\nu}}{(1-q)^{\nu}}\sum_{n=0}^{\infty}\frac{q^{n^2+n\nu}(y/2)^{2n}}{(1-q)^n(q;q)_n\Gamma_q(\nu+n+1)},~~~y\in\C.
$$
The sequences $g^{(1)}_n=(q;q)_n^{-1}(1-q)^{-n}$ and $g^{(2)}_n=q^{n^2+n\nu}(q;q)_n^{-1}(1-q)^{-n}$, $n=0,1,\ldots$, are immediately seen to be doubly positive, and hence we are in the position to apply Theorem~\ref{th:gammadenom1} and its corollaries to the functions $(y/2)^{-\nu}(1-q)^{\nu}I_{\nu}^{(j)}(y)$, $j=1,2$, with $x=(y/2)^2$ and $\mu=\nu+1$.
In particular, Corollary~\ref{cr:Delta-g-bounds} yields the following bounds:
\begin{multline*}
\frac{(y/2)^{2\nu+\alpha+\beta}[(q^{\nu+1};q)_{\alpha}^{-1}-(q^{\nu+\beta+1};q)_{\alpha}^{-1}]}{\Gamma_q(\nu+1)\Gamma_q(\nu+\beta+1)(1-q)^{2\nu+\beta}}
\le I^{(j)}_{\nu+\alpha}(y;q)I^{(j)}_{\nu+\beta}(y;q)-I^{(j)}_{\nu}(y;q)I^{(j)}_{\nu+\alpha+\beta}(y;q)
\\
\leq\left[1-\frac{(q^{\nu+1};q)_{\alpha}}{(q^{\nu+\beta+1};q)_{\alpha}}\right]
I^{(j)}_{\nu+\alpha}(y;q)I^{(j)}_{\nu+\beta}(y;q)
\end{multline*}
for $\nu>\max(-\beta-1,-2)$,  $\beta>0$, $\alpha\in\N$ and $0<y<R_j$, where $R_1=2$, $R_2=\infty$.  For $\alpha=\beta=1$ and changing $\nu$ to $\nu-1$ this leads to the Tur\'{a}n type inequality (recall that $[x]_q=(1-q^x)/(1-q)$),
$$
\frac{(y/2)^{2\nu}q^{\nu}(1-q)^{-2\nu}}{[\nu+1]_q[\Gamma_q(\nu+1)]^2}\leq (I^{(j)}_{\nu}(y;q))^2-I^{(j)}_{\nu-1}(y;q)I^{(j)}_{\nu+1}(y;q)\leq\frac{q^{\nu}}{[\nu+1]_q}(I^{(j)}_{\nu}(y;q))^2.
$$
Using the limit relations \cite[Remark~3.1]{OR1} $\lim_{q\uparrow1}I^{(j)}_{\nu}((1-q)y;q)=I_{\nu}(y)$, $j=1,2$, the above inequality reduces to \cite[(26)]{KK1}.  Finally, we remark that, in fact, $\nu\to{I^{(1)}_{\nu}(y;q)}$ is log-concave on $(-1,\infty)$ for each $0<y<2$. A proof of this claim will be presented in our forthcoming work.

\medskip

\textbf{Example~2.} The $q$-Kummer function can either be defined by $r=s=1$ case of (\ref{eq:r-phi-s}):
$$
\phi^{(1)}(a;b;z)={}_1\phi_1(a;b;q;-z)=\sum_{n=0}^{\infty}\frac{(a;q)_nq^{n(n-1)/2}}{(b,q;q)_n}z^n,~~~z\in\C,
$$
which corresponds to the Gasper-Rahman definition \cite[(1.2.22)]{GR} or by the series
$$
\phi^{(2)}(a;b;z)={}_2\phi_1(a,0;b;q;z)=\sum_{n=0}^{\infty}\frac{(a;q)_n}{(b,q;q)_n}z^n,~~~|z|<1,
$$
which is the definition used by Bailey and Slater \cite[(3.2.1.11)]{Slater}.  Both sequences $f_n^{(1)}=q^{n(n-1)/2}(b;q)_n^{-1}$ and $f_n^{(2)}=(b;q)_n^{-1}$ are easily seen to be doubly positive for any $0<b<1$.  Hence, we are in the position to apply Theorem~\ref{th:pochhammertop} and its corollaries and conclude that $\mu\to\phi^{(j)}(q^{\mu};b;q;z)$ is log-concave on $[0,\infty)$ for $j=1,2$ and the power series coefficients of the generalized Tur\'{a}nians $\Delta_{\phi^{(j)}}(\alpha,\beta;x)$ (see definition (\ref{eq:delta-defined})) have positive coefficients at all positive powers of $x$. Furthermore, by Corollary~\ref{cr:Delta-f-bounds} we have the following bounds:
$$
\frac{(1-q^{\alpha})(1-q^{\beta})}{(1-b)(1-q)}q^{\mu}x\!\le\!\Delta_{\phi^{(j)}}(\alpha,\beta;x)\!<\!\left[1-\frac{\Gamma_q(\mu+\alpha)\Gamma_q(\mu+\beta)}{\Gamma_q(\mu)\Gamma_q(\mu+\alpha+\beta)}\right]
\!\phi^{(j)}(q^{\mu+\alpha};b;x)\phi^{(j)}(q^{\mu+\beta};b;x)
$$
for $j=1,2$.  These results for the function $\phi^{(2)}(a;b;z)$ strengthen and generalize \cite[Theorem~3.2]{BRS} due to
Baricz, Raghavendar and Swaminathan which asserts that
$$
\phi^{(2)}(q^a,q^c;x)^2>\phi^{(2)}(q^{a+m},q^{c};x)\phi^{(2)}(q^{a-m},q^{c};x)
$$
for $x>0$, $c>0$, $a\geq{m-1}$, $m\in\N$.

\medskip

\textbf{Example~3.}  Keeping the notation of Example~2, we see that the functions $h_j(b;x)=\phi^{(j)}(a;b;x)$, $j=1,2$, satisfy the conditions of Theorem~\ref{th:q-poch-denom} and its corollaries.  In particular, $\mu\to\phi^{(j)}(a;q^{\mu};x)$ is log-convex on $(0,\infty)$, $\Delta_{h_j}(\alpha,\beta;x)$ has negative power series coefficients
for all $\mu,\alpha,\beta>0$.  Furthermore, it is easy to check that the sequences $h_n^{(1)}=q^{n(n-1)/2}(a;q)_n/(q;q)_n$ and  $h_n^{(2)}=(a;q)_n/(q;q)_n$, $n\in\N_0$, are both doubly positive if $0<a<q$. Hence, by Corollary~\ref{cr:Delta-h-bounds} we have for $j=1,2$:
\begin{multline*}
\left[\frac{\Gamma_q(\mu+\alpha)\Gamma_q(\mu+\beta)}{\Gamma_q(\mu)\Gamma_q(\mu+\alpha+\beta)}-1\right]\phi^{(j)}(a;q^{\mu};x)\phi^{(j)}(a;q^{\mu+\alpha+\beta};x)<\Delta_{h_j}(\alpha,\beta;x)
\\
\le\frac{-(1-a)xq^{\mu}(1-q^{\alpha})(1-q^{\beta})}{(1-q)(1-q^{\mu})(1-q^{\mu+\alpha+\beta})(1-q^{\mu+\alpha})(1-q^{\mu+\beta})}
\end{multline*}
for all $0<a<q$, $\mu,\beta>0$, $\alpha\in\N$ and $0\le{x}<R_j$, where $R_1=\infty$, $R_2=1$. The upper bound holds all $0<a<1$ and $\mu,\alpha,\beta>0$.
Note that for $j=2$ this inequality is a refinement of the second statement in \cite[Theorem~1]{MehrSitn}. It is interesting to remark here that the first statement in \cite[Theorem~1]{MehrSitn} states essentially that the ratio in the middle of the first inequality in Corollary~\ref{cr:Delta-h-bounds} built on $h_2(b;x)=\phi^{(2)}(a;b;x)$ is a monotone function of $x$.

\medskip

\textbf{Example~4.} Put $f_{n}=(b;q)_n/(c;q)_n$ for $n=0,1,\ldots$.  Then for $f(a;x)$ defined in (\ref{eq:f}) we have
$$
f(a;x)={_2}\phi_1(a,b;c;q;x),~~~|x|<1,
$$
where the function on the right hand side is given in (\ref{eq:r-phi-s}).  It is straightforward to check that $\{f_n\}_{n\ge0}$ is doubly positive if $0<b<c<1$. Under this condition we can apply Theorem~\ref{th:pochhammertop} and its corollaries.  In particular, the function $\mu\to{_2}\phi_1(q^{\mu},b;c;q;x)$ is log-concave on $[0,\infty)$ for each fixed $0<x<1$ and $\Delta_f(\alpha,\beta;x)$ has positive power series coefficients. Furthermore, according to Corollary~\ref{cr:Delta-f-bounds} we have:
$$
\frac{q^{\mu}x(1-b)(1-q^{\alpha})(1-q^{\beta})}{(1-q)(1-c)}\!\le\!\Delta_f(\alpha,\beta;x)\!<\!\left[1-\frac{\Gamma_q(\mu+\alpha)\Gamma_q(\mu+\beta)}{\Gamma_q(\mu)\Gamma_q(\mu+\alpha+\beta)}\right]
\!f(q^{\mu+\alpha};x)f(q^{\mu+\beta};x)
$$
for all $\alpha,\beta>0$, $\mu\ge0$ and $0\le{x}<1$.

\medskip

\textbf{Example~5.} Put $h_{n}=(a;q)_n(b;q)_n$ for $n\in\N_0$.  Then for $h(c;x)$ defined in (\ref{eq:h}) we have
$$
h(c;x)={_2}\phi_1(a,b;c;q;x),~~~|x|<1,
$$
where the function on the right hand side is given in (\ref{eq:r-phi-s}).  Hence, we are in the position to apply Theorem~\ref{th:q-poch-denom} to conclude that
the function $\mu\to{_2}\phi_1(a,b;q^{\mu};q;x)$ is log-convex on $[0,\infty)$ for any fixed $0<a,b,x<1$  and $\Delta_h(\alpha,\beta;x)$ defined in (\ref{eq:Delta-h}) has negative power series coefficients at all positive powers of $x$.  Moreover, according to the upper bound in the second inequality of Corollary~\ref{cr:Delta-h-bounds} we have
$$
\Delta_h(\alpha,\beta;x)\le\frac{-(1-a)(1-b)(1-q^{\alpha})(1-q^{\beta})q^{\mu}x}{(1-q^{\mu})(1-q^{\mu+\alpha+\beta})(1-q^{\mu+\alpha})(1-q^{\mu+\beta})}.
$$

\medskip

Before turning to the next example let us define the rational function
\begin{equation}\label{eq:Rpq}
R_{r,s}(x)=\frac{\prod_{k=1}^{r}(a_k+x)}{\prod_{k=1}^{s}(b_k+x)}
\end{equation}
with  positive $a_k,b_k$. Let $e_m(\mathbf{c})=e_m(c_1,\dots,c_r)$ denote the $m$-th elementary symmetric polynomial,
$$
e_0(c_1,\dots,c_r)=1,~~~~~e_1(c_1,\dots,c_r)=c_1+c_2+\dots+c_r,
$$
$$
e_2(c_1,\dots,c_r)=c_1c_2+c_1c_3+\dots+c_1c_q+c_2c_3+\dots+c_2c_q+\dots+c_{r-1}c_{r},\ldots,
$$
$$
e_r(c_1,\dots,c_r)=c_1c_2\cdots{c_r}.
$$
We will need the following lemma  \cite[Lemma~3]{KK4}.

\begin{lemma}\label{lem:incrdecr}
If $r\le{s}$ and
\begin{equation}\label{eq:decr}
\frac{e_s(\b)}{e_r(\a)}\le\frac{e_{s-1}(\b)}{e_{r-1}(\a)}\le\dots\le\frac{e_{s-r+1}(\b)}{e_1(\a)}\leq e_{s-r}(\b),
\end{equation}
then the function $R_{r,s}(x)$ defined in \emph{(\ref{eq:Rpq})} is monotone decreasing on $(0,\infty)$.
\end{lemma}

\medskip

\textbf{Example 6.}  Suppose $\alpha_i>0$, $i=1,\ldots,r$,  $\beta_j>0$, $j=1,\ldots,s$ and $r\le{s}$. Define $f_{n}=q^{n(n-1)(s-r)/2}(q^{\alpha_1},\dots,q^{\alpha_r};q)_n/(q^{\beta_1},\dots,q^{\beta_s};q)_n$ for $n\in\N_0$.
This sequence is clearly positive.  By definition of log-concavity, it is doubly positive if and only if the sequence $\{f_{n+1}/f_{n}\}_{n\ge0}$ is decreasing.
Straightforward calculation then reveals:
\begin{multline*}
\frac{f_{n+1}}{f_{n}}\!=\!q^{(s-r)n}\frac{(1-q^{\alpha_1+n})\cdots(1-q^{\alpha_r+n})}{(1-q^{\beta_1+n})\cdots(1-q^{\beta_s+n})}\!
\\
=\!q^{(s-r)n}\frac{\prod_{k=1}^{r}q^{\alpha_k}}{\prod_{k=1}^{s}q^{\beta_k}}\frac{\prod_{k=1}^{r}(q^{-\alpha_k}-q^n)}{\prod_{k=1}^{s}(q^{-\beta_k}-q^n)}
\!=\!q^{(s-r)n}\frac{\prod_{k=1}^{r}q^{\alpha_k}}{\prod_{k=1}^{s}q^{\beta_k}}\frac{\prod_{k=1}^{r}(a_k+y_n)}{\prod_{k=1}^{s}(b_k+y_n)},
\end{multline*}
where $a_k=q^{-\alpha_k}-1>0$, $b_k=q^{-\beta_k}-1>0$, $y_n=1-q^n>0$. As $n\to{y_n}$ is increasing, it is clear that $\{f_{n+1}/f_{n}\}_{n\ge0}$ is decreasing if
$$
F(y)=\frac{\prod_{k=1}^{r}(a_k+y)}{\prod_{k=1}^{s}(b_k+y)}
$$
is decreasing.  Hence, if condition (\ref{eq:decr}) is satisfied with $a_k=q^{-\alpha_k}-1$, $b_k=q^{-\beta_k}-1$, then
$$
f(a;x)=\sum\limits_{n=0}^{\infty}f_n\frac{(a;q)_n}{(q;q)_n}x^n={}_{r+1}\phi_{s}\left(a,q^{\alpha_1},\ldots,q^{\alpha_r};q^{\beta_1},\ldots,q^{\beta_s};q;(-1)^{s-r}x\right)
$$
satisfies Theorem~\ref{th:pochhammertop} and Corollaries~\ref{cr:Delta-fCM} and \ref{cr:Delta-f-bounds}.  Note that the logarithmic convexity of ${}_{r+1}\phi_{r}$ with respect to shifts in $\beta_k$, $k\in\{1,\ldots,r\}$, has been recently proved in \cite[Theorem~2]{MehrSitn}.

\medskip

\textbf{Example 7.}  Suppose $\alpha_i>0$, $i=1,\ldots,r$,  $\beta_j>0$, $j=1,\ldots,s-1$, $\beta_s=1$ and $r\le{s+1}$. Define $g_{n}=q^{n(n-1)(1+s-r)/2}(q^{\alpha_1},\dots,q^{\alpha_r};q)_n/(q^{\beta_1},\dots,q^{\beta_s};q)_n$ for $n\in\N_0$.
Similarly to the previous example this sequence is doubly positive if (\ref{eq:decr}) is satisfied with $a_k=q^{-\alpha_k}-1$, $k=1,\ldots,r$ and $b_k=q^{-\beta_k}-1$, $k=1,\ldots,s$. Then using (\ref{eq:qGamma-prop}) we have
$$
g(\mu;x)=\sum\limits_{n=0}^{\infty}\frac{g_nx^n}{\Gamma_q(\mu+n)}=\frac{1}{\Gamma_q(\mu)}
{}_{r}\phi_{s}\left(q^{\alpha_1},\ldots,q^{\alpha_r};q^{\mu},q^{\beta_1},\ldots,q^{\beta_{s-1}};q;(1-q)(-1)^{1+s-r}x\right)
$$
satisfies Theorem~\ref{th:gammadenom1} and Corollaries~\ref{cr:g-multconc}, \ref{cr:cm2} and \ref{cr:Delta-g-bounds}.

\paragraph{Acknowledgements} This research has been supported by the Russian Science Foundation under project 14-11-0002.

\end{document}